\documentclass[11pt, hidelinks,letterpaper, oneside, reqno]{amsart}
\usepackage{amsfonts, amsmath, amsthm, amssymb}
\usepackage{mathtools}
\usepackage{mathrsfs}
\usepackage{tikz-cd}
\usepackage{hyperref}
\usepackage[dvipsnames]{xcolor}
\usepackage[utf8]{inputenc}
\usepackage{newpxtext}
\usepackage[T1]{fontenc}

\hypersetup{
    colorlinks=true,
    linkcolor=Blue,
    filecolor=Blue,
    urlcolor=Blue,
    citecolor=Blue,
}
\usepackage[top=2cm, bottom=2cm, left=2cm, right=2cm]{geometry}

\newcommand{\too}{\longrightarrow}
\newcommand{\defn}[1]{\emph{\textbf{#1}}}

\DeclareMathOperator{\id}{id}
\makeatletter
\newtheoremstyle{plain}%
  {2pt}{2pt}%
  {\itshape}{0pt}{\bfseries}{.}{0.5em}%
  {\thmname{#1}\thmnumber{ #2}\thmnote{ \normalfont\itshape(#3)}}
\newtheoremstyle{definition}%
  {2pt}{2pt}%
  {\upshape}{0pt}{\bfseries}{.}{0.5em}%
  {\thmname{#1}\thmnumber{ #2}\thmnote{ \normalfont\itshape(#3)}}
\newtheoremstyle{remark}%
  {2pt}{2pt}%
  {\upshape}{0pt}{\bfseries}{.}{0.5em}%
  {\thmname{#1}\thmnumber{ #2}\thmnote{ \normalfont\itshape(#3)}}
\makeatother

\theoremstyle{plain}
\newtheorem{theorem}{Theorem}[section]
\newtheorem{lemma}[theorem]{Lemma}
\newtheorem{proposition}[theorem]{Proposition}
\newtheorem{corollary}[theorem]{Corollary}

\theoremstyle{definition}
\newtheorem{definition}[theorem]{Definition}
\newtheorem{example}[theorem]{Example}

\theoremstyle{remark}
\newtheorem{remark}[theorem]{Remark}

\makeatletter
\renewcommand\section{\@startsection{section}{1}%
  \z@{.7\linespacing\@plus\linespacing}{.5\linespacing}%
  {\normalfont\bfseries\large\centering}}
\renewcommand\subsection{\@startsection{subsection}{2}%
  \z@{.5\linespacing\@plus.7\linespacing}{-.5em}%
  {\normalfont\bfseries\centering}}
\makeatother

\usepackage[style=numeric, sorting=nyt, backend=biber]{biblatex}
\address{Yassine Ait Mohamed}

\begin{document}

\author{Yassine Ait Mohamed}
\title[Generalized Homogeneous Derivations on Graded Rings]{Generalized Homogeneous Derivations on Graded Rings}

\begin{abstract}
We introduce a notion of generalized homogeneous derivations on graded rings, extending the homogeneous derivations of Kanunnikov, and single out the gr-generalized derivations among them, those preserving the degree of every homogeneous component. Several results originally established for prime rings are extended to gr-prime rings. We also study the algebraic and module-theoretic structure of these maps, establish their functorial properties, and construct categories that record their derivation structure in both ring and module contexts.
\end{abstract}



\maketitle

\thispagestyle{empty}

\section{Introduction}\label{LV0PRKG}

The concept of generalized derivations, first introduced by Bre\v{s}ar \cite{Br91}, has been widely studied since the work of Hvala \cite{Hava}. Many results about derivations have been extended to the context of generalized derivations (see, for example, \cite{dhara2025action,HuangDavvaz13,LeeSh,ZhaoLiu2024}). In 2018, Kanunnikov \cite{GradQ} introduced the concept of homogeneous derivations for graded rings. Homogeneous derivations are derivations in the classical sense that are compatible with the graded structure of the ring. The purpose of this paper is to introduce a novel concept called generalized homogeneous derivations for graded rings, which extends the notion of homogeneous derivations to the context of generalized derivations while respecting the graded structure. We establish that homogeneous derivations are a special case of generalized homogeneous derivations, thereby providing a unifying framework for the study of derivations in graded rings. As a consequence of this new concept, we also extend several important results, originally proven for prime rings, to the setting of gr-prime rings.

The present paper continues our earlier work \cite{yassine}, in which several classical results on derivations were extended to the graded setting.

Let $\Gamma$ denote a group with identity element $e$. A ring $R$ is called $\Gamma$-graded if it can be decomposed as $R = \bigoplus_{\tau\in \Gamma} R_\tau$ into additive subgroups such that $R_{\tau} R_{\sigma}\subset R_{\tau\sigma}$ for all $\tau, \sigma \in \Gamma$. The collection of homogeneous elements $\mathcal{H}(R) = \bigcup_{\tau\in \Gamma} R_\tau$ consists of elements $a \in R_\tau$ having degree $\deg a = \tau$. Each element $x \in R$ has a unique representation $x = \sum_{\tau\in \Gamma} x_\tau$ where $x_\tau \in R_\tau$ are the homogeneous components. The graded structure naturally extends to tensor products through
\[
(R \otimes_K S)_{\gamma} = \bigoplus_{\tau\sigma = \gamma} R_{\tau} \otimes_K S_{\sigma}
\]
for $\Gamma$-graded $K$-algebras and similarly to polynomial rings.

An ideal $\mathfrak{a} \subseteq R$ is graded when $\mathfrak{a} = \bigoplus_{\tau\in \Gamma} (\mathfrak{a} \cap R_\tau)$, where we denote $\mathfrak{a}_\tau = \mathfrak{a} \cap R_\tau$. Ring homomorphisms $\varphi: R \too S$ between $\Gamma$-graded rings are graded if $\varphi(R_\tau) \subseteq S_\tau$. We write $\operatorname{Hom}(R,S)^{gr}$ for the set of all graded homomorphisms. When $\mathfrak{a}$ is graded, quotient rings inherit the canonical grading structure $(R/\mathfrak{a})_{\tau} := \{\bar{r} \in R/\mathfrak{a} : r \in R_{\tau}\}$. A graded ring $R$ is gr-prime if $aRb = 0$ implies $a = 0$ or $b = 0$ for homogeneous elements $a, b \in \mathcal{H}(R)$, and gr-semiprime if $aRa = 0$ implies $a = 0$ for $a \in \mathcal{H}(R)$.

Given $\{R_i\}_{i\in I}$ a finite collection of $\Gamma_i$-graded rings, set $R=\prod_{i\in I}R_i$. For $(\tau_i)_i\in\prod_{i\in I}\Gamma_i$, put $R_{(\tau_i)_i} := \prod_{i\in I}(R_i)_{\tau_i}$. Then
\begin{equation}\label{PRODGRAD}
R = \bigoplus_{(\tau_i)_i\,\in\,\prod_{i\in I}\Gamma_i} R_{(\tau_i)_i}
\end{equation}
is a $\prod_{i\in I}\Gamma_i$-graded ring.

A $\Gamma$-graded module over a $\Gamma$-graded ring $R$ decomposes as $M = \bigoplus_{\tau \in \Gamma} M_\tau$ with the compatibility condition $R_\sigma \cdot M_\tau \subseteq M_{\sigma\tau}$. Graded homomorphisms satisfy $f(M_\tau) \subseteq N_\tau$, and tensor products exhibit the multiplicative grading
\[
(M \otimes_R N)_k = \bigoplus_{\tau\sigma = k} M_\tau \otimes_R N_\sigma.
\]

An additive mapping $d: R \too R$ is a derivation if it satisfies the Leibniz rule $d(xy)=d(x)y+xd(y)$ for all $x,y \in R$. A derivation $d$ is homogeneous if $d(\mathcal{H}(R))\subseteq \mathcal{H}(R)$ \cite{GradQ}. Inner derivations have the form $d(x) = [a, x]$ for some fixed $a \in R$, where $[a,x] = ax - xa$ denotes the commutator. A generalized derivation is an additive mapping $F : R\too R$ satisfying $F(xy) = F(x)y + xd(y)$ for all $x,y \in R$, where $d$ is the associated derivation of $F$.

\medskip
Throughout this paper, we adopt the following conventions.
\begin{itemize}
    \item All polynomial rings $\mathbb{C}[t_1, \ldots, t_n]$ are equipped with the standard $\mathbb{Z}$-grading by total degree, where $\deg(t_1^{a_1} \cdots t_n^{a_n}) = a_1 + \cdots + a_n$.
    \item \emph{The $\pm$ notation:} When a condition involves the symbol $\pm$, such as
    \[
    F(xy) \pm xy \in Z(R) \quad \text{or} \quad F_1(x)F_2(y) \pm xy \in Z(R),
    \]
    we mean that at least one of the two possibilities holds. The case with $-$ always reduces to the case with $+$ by replacing $F$ (or $F_1, F_2$) by $-F$ (or $-F_1, -F_2$), so the proofs establish the result for one sign and invoke this reduction for the other.
    \item \emph{Abelian grading group:} We restrict our attention to abelian grading groups $\Gamma$ throughout. This is necessary: in general the Lie bracket $[x,y] = xy - yx$ and the Jordan product $x \circ y = xy + yx$ do not preserve homogeneity, as Example~\ref{PDPFF5E} shows. Abelianness of $\Gamma$ is what guarantees preservation of homogeneity, and the commutator-based arguments of this paper rely on it throughout.
\end{itemize}

\begin{example}\label{PDPFF5E}
   Let $R=M_{4}(k)$ denote the ring of $4\times 4$ matrices over a field $k$, and let $D_{10}=\langle a,b\mid a^{5} =b^{2}=e, \, bab=a^{-1}\rangle$ be the dihedral group of order $10$. We define a $D_{10}$-grading on $R$ by setting
   \[
   R_{e}:=\begin{pmatrix} k&0&0&0 \\ 0&k&0&0 \\ 0&0&k&0 \\ 0&0&0&k \end{pmatrix}, \quad
   R_{a}:=\begin{pmatrix} 0&k&0&0 \\ 0&0&k&0 \\ 0&0&0&0 \\ 0&0&0&0 \end{pmatrix}, \quad
   R_{a^{2}}:=\begin{pmatrix} 0&0&k&0 \\ 0&0&0&0 \\ 0&0&0&0 \\ 0&0&0&0 \end{pmatrix},
   \]
   \[
   R_{a^{3}}:=\begin{pmatrix} 0&0&0&0 \\ 0&0&0&0 \\ k&0&0&0 \\ 0&0&0&0 \end{pmatrix}, \quad
   R_{b}:=\begin{pmatrix} 0&0&0&0 \\ 0&0&0&k \\ 0&0&0&0 \\ 0&k&0&0 \end{pmatrix}, \quad
   R_{ab}:=\begin{pmatrix} 0&0&0&k \\ 0&0&0&0 \\ 0&0&0&0 \\ k&0&0&0 \end{pmatrix},
   \]
   \[
   R_{a^{4}b}:=\begin{pmatrix} 0&0&0&0 \\ 0&0&0&0 \\ 0&0&0&k \\ 0&0&k&0 \end{pmatrix}, \quad
   R_{a^{4}}:=\begin{pmatrix} 0&0&0&0 \\ k&0&0&0 \\ 0&k&0&0 \\ 0&0&0&0 \end{pmatrix}, \quad
   R_{a^{2}b}=R_{a^{3}b}=0.
   \]
   This is the elementary grading determined by the sequence $(\sigma_1,\sigma_2,\sigma_3,\sigma_4)=(e,a,a^2,ab)$ attached to the four rows/columns, via $\deg(e_{ij})=\sigma_i^{-1}\sigma_j$; one checks directly that this assignment reproduces exactly the components listed above. One verifies by direct computation that, for $x=e_{24}+e_{42}\in R_b$ and $y=e_{21}+e_{32}\in R_{a^4}$,
   \[
   [x,y]=e_{41}-e_{34}, \quad x\circ y = e_{41}+e_{34},
   \]
   and that, as $x$ and $y$ range over $R_b$ and $R_{a^4}$, both $[x,y]$ and $x\circ y$ span the same two-dimensional space
   \[
   \begin{pmatrix} 0&0&0&0 \\ 0&0&0&0 \\ 0&0&0&k \\ k&0&0&0 \end{pmatrix},
   \]
   which is not contained in $\mathcal{H}(R)$, since a generic element of it, such as $e_{34}+e_{41}$, has nonzero components in the two distinct degrees $a^4b$ and $ab$.
\end{example}

The paper is organized as follows. Section \ref{CWNPGW9} introduces generalized homogeneous derivations and isolates the gr-generalized derivations among them. Section \ref{FUNCTPROP} develops their algebraic and categorical structure. Section \ref{9J1ICTX} establishes commutativity criteria for gr-prime rings. Section~\ref{YC4OM3W} extends the theory to graded modules.
\section{Generalized homogeneous derivations}\label{CWNPGW9}

We introduce generalized homogeneous derivations, show that the naive definition fails to carry an additive or product structure, and isolate the gr-generalized derivations that repair both failures.

\subsection{Definition of generalized homogeneous derivations}\label{SUB-DEFGHD}

Throughout this subsection, $R$ denotes a ring graded by a group $\Gamma$, not assumed abelian unless stated otherwise.

\begin{definition}\label{DD467KD}
An additive mapping $F : R \too R$ is called a \defn{generalized homogeneous derivation} if there exists a homogeneous derivation $d : R \too R$ such that
\begin{enumerate}
    \item[(1)] $F(xy) = F(x)y + x d(y)$ for all $x,y \in R$;
    \item[(2)] $F(r) \in \mathcal{H}(R)$ for all $r \in \mathcal{H}(R)$.
\end{enumerate}
The mapping $d$ is called an associated homogeneous derivation of $F$. We denote such a generalized homogeneous derivation by $(F, d)_h$, where the subscript `$h$' records the homogeneity condition, and write $\mathfrak{Der}^{gh}_\Gamma(R)$ for the collection of all generalized homogeneous derivations of $R$.
\end{definition}

\begin{example}\label{EX-MNCT}
Let $R = M_n(\mathbb{C}[t])$ with the $\mathbb{Z}_2$-grading where $R_0$ consists of matrices with polynomial entries having only even-degree monomials, and $R_1$ consists of matrices with polynomial entries having only odd-degree monomials. Define $d : R \too R$ by $d(A) = \dfrac{\partial A}{\partial t}$ (entrywise differentiation) and $F : R \too R$ by $F(A) = tA + d(A)$. Then $(F,d)_h$ is a generalized homogeneous derivation.
\end{example}

\begin{proposition}\label{PROP-STRICT}
Let $R$ be a nontrivially $\Gamma$-graded ring. Then the following inclusions hold:
\[
\mathfrak{Der}^{h}_\Gamma(R) \subsetneq \mathfrak{Der}^{gh}_\Gamma(R) \subsetneq \mathfrak{Gen}(R),
\]
where $\mathfrak{Der}^{h}_\Gamma(R)$ and $\mathfrak{Gen}(R)$ denote the sets of homogeneous derivations and generalized derivations on $R$, respectively. Both inclusions are strict.
\end{proposition}

\begin{proof}
The inclusions follow directly from the definitions. For strictness, consider $R = \mathbb{C}[t_1,t_2]$ with the standard $\mathbb{Z}$-grading. Define
\[
F(f) = t_1 f + t_1 t_2 \frac{\partial f}{\partial t_1}
\quad\text{and}\quad
d(f) = t_1 t_2 \frac{\partial f}{\partial t_1}.
\]
Then $(F,d)_h \in \mathfrak{Der}^{gh}_\Gamma(R)$, but $F \notin \mathfrak{Der}^{h}_\Gamma(R)$, since $F$ itself is not a derivation. For the second inclusion, define $G(f) = t_1 f + \frac{\partial f}{\partial t_2}$. Then $G \in \mathfrak{Gen}(R)$, with associated derivation $\frac{\partial}{\partial t_2}$, but $G \notin \mathfrak{Der}^{gh}_\Gamma(R)$, since $G(1) = t_1$ is homogeneous while $G(t_2) = t_1t_2 + 1$ is not.
\end{proof}

\begin{proposition}\label{T3UEA3G}
Let $R$ be a $\Gamma$-graded ring. Then $R$ admits a nonzero generalized homogeneous derivation if any of the following holds:
\begin{enumerate}
\item $R$ has a nonzero homogeneous derivation;
\item $R_{\sigma} \cap Z(R) \neq 0$ for some $\sigma \in \Gamma$;
\item $C_R(R_e) = R$ and $R_e \neq 0$.
\end{enumerate}
\end{proposition}

\begin{proof}
(1) If $d \neq 0$ is a homogeneous derivation, then $(d,d)_h$ is a nonzero generalized homogeneous derivation by definition.

(2) For a nonzero $a \in R_{\sigma} \cap Z(R)$, define $F_a(r) = ar$. Since $a$ is central and homogeneous, $F_a$ preserves homogeneous elements, and $(F_a, 0)_h$ satisfies the required conditions.

(3) For a nonzero $b \in R_e$ with $C_R(R_e) = R$, define $F_b(r) = br$. Since $b$ is central, $(F_b, 0)_h$ is well defined and nonzero.
\end{proof}

Two familiar facts about generalized derivations fail once homogeneity is imposed. For ordinary generalized derivations $F_1,F_2$ of $R$, the sum $F_1+F_2$ is again a generalized derivation, and a finite family of rings, each carrying a generalized derivation, endows the product ring with a generalized derivation defined coordinatewise. Neither property survives for generalized \emph{homogeneous} derivations once the grading is nontrivial, as the next two examples show.

\begin{example}\label{EVB9OOA}
Consider the polynomial ring $\mathbb{C}[t_1, t_2, t_3]$ equipped with the standard $\mathbb{Z}$-grading, and define
\[
F_1(f) = d_1(f) = t_3 \frac{\partial f}{\partial t_1}
\quad\text{and}\quad
F_2(f) = d_2(f) = \frac{\partial f}{\partial t_2}.
\]
Both $(F_1,d_1)_h$ and $(F_2,d_2)_h$ lie in $\mathfrak{Der}^{gh}_{\mathbb Z}(R)$, yet
\[
(F_1 + F_2)(f) = t_3 \frac{\partial f}{\partial t_1} + \frac{\partial f}{\partial t_2}
\]
sends the homogeneous element $t_1t_2$ to $t_2t_3+t_1$, which is not homogeneous. So $\mathfrak{Der}^{gh}_\Gamma(R)$ carries no natural additive structure.
\end{example}

\begin{example}\label{PRODFAIL-EX}
For the second property, take $R_1=R_2=\mathbb C[t]$ with the standard $\mathbb Z$-grading, and set $F_1(f)=tf$, $d_1=0$, and $F_2(f)=f$, $d_2=0$. Both $(F_1,0)_h$ and $(F_2,0)_h$ lie in $\mathfrak{Der}^{gh}_{\mathbb Z}(\mathbb C[t])$. For the natural componentwise grading $(R_1\times R_2)_n=(R_1)_n\times(R_2)_n$, the element $(1,1)$ is homogeneous of degree $0$, yet $F(1,1)=(t,1)$ lies in no homogeneous component of $R_1\times R_2$.
\end{example}

\subsection{Definition of gr-generalized derivations}\label{SUB-DEFGRGEN}

The failure recorded in Examples \ref{EVB9OOA} and \ref{PRODFAIL-EX} comes from the same source: Definition \ref{DD467KD} only asks that $F$ and $d$ send $\mathcal H(R)$ into itself, not that they preserve the degree of each homogeneous component. Restricting to maps that do preserve degree repairs both properties, and singles out the class of derivations we study for the rest of the paper.

\begin{definition}\label{37Q9AH8}
Let $R$ be a ring graded by an arbitrary group $\Gamma$. A generalized homogeneous derivation $(F,d)_h$ is called a \defn{gr-generalized derivation} if
\[
F(R_\tau) \subseteq R_\tau \quad \text{and} \quad d(R_\tau) \subseteq R_\tau \quad \text{for all } \tau \in \Gamma.
\]
The set of all such derivations is denoted by $p\mathfrak{Der}^{gh}_\Gamma(R)$.
\end{definition}

\begin{example}\label{EX-GRGEN-EULER}
Let $R = \mathbb{C}[t_1,t_2]$ be the polynomial ring with the standard $\mathbb{Z}$-grading. Define $F(f) = d(f) = t_1 \frac{\partial f}{\partial t_1} + t_2 \frac{\partial f}{\partial t_2}$ for all $f \in R$. Then $(F,d)_h$ is a gr-generalized derivation, since $F$ acts on a monomial of total degree $n$ by multiplication by $n$.
\end{example}

The hierarchy of notions considered in this paper is as follows:

\[
\begin{tikzcd}[ampersand replacement=\&,cramped,row sep=scriptsize]
	\&\& {\text{derivations}} \&\& \\
	{\text{homogeneous derivations}} \&\&\&\& {\text{generalized derivations}} \\
	\&\& {\text{generalized homogeneous derivations}} \\
	\\
	\&\& {\text{gr-generalized derivations}}
	\arrow[Rightarrow, from=1-3, to=2-5]
	\arrow[Rightarrow, from=2-1, to=1-3]
	\arrow[Rightarrow, from=2-1, to=3-3]
	\arrow[Rightarrow, from=3-3, to=2-5]
	\arrow[Rightarrow, from=5-3, to=3-3]
\end{tikzcd}
\]

\begin{remark}
Every gr-generalized derivation of $R$ restricts to a generalized derivation of the identity component $R_e$; the converse fails, as the next example shows.
\end{remark}

\begin{example}
Let $A=\mathbb{C}[t]$ and $R=M_2(A)$, graded by $R_0=\{\mathrm{diag}(f,g)\}$, $R_1=\{fe_{12}+ge_{21}\}$. Identify $R_0\cong A\times A$ via $\mathrm{diag}(f,g)\leftrightarrow(f,g)$ and set $d_e(f,g)=(\frac{\partial f }{\partial t} ,0)$. Then $(d_e,d_e)$ is a generalized derivation of $R_0=R_e$.

Suppose some $(F,d)_h\in p\mathfrak{Der}^{gh}_{\mathbb{Z}_2}(R)$ restricts to $d_e$ on $R_0$. Since $d(R_1)\subseteq R_1$, write $d(e_{12})=pe_{12}+qe_{21}$, $p,q\in A$. As $d(e_{11})=d_e(1,0)=0$, Leibniz on $e_{12}=e_{11}e_{12}$ gives $d(e_{12})=e_{11}d(e_{12})=pe_{12}$, so $q=0$. Since $d(te_{11})=e_{11}$ and $d(te_{22})=0$, expanding $te_{12}=(te_{11})e_{12}=e_{12}(te_{22})$ via Leibniz yields $(1+tp)e_{12}=tp\,e_{12}$, i.e., $e_{12}=0$, a contradiction. Hence no gr-generalized derivation of $R$ restricts to $d_e$ on $R_e$.
\end{example}

\section{Functorial and algebraic properties of gr-generalized derivations}\label{FUNCTPROP}

We now examine how generalized homogeneous derivations transport along graded ring homomorphisms, the algebraic structure carried by $p\mathfrak{Der}^{gh}_\Gamma(R)$, and the behavior of gr-generalized derivations under finite products, polynomial extension, and tensor product.

\begin{definition}\label{GRDIFFIDEAL}
Let $(F,d)_h$ be a generalized homogeneous derivation of a $\Gamma$-graded ring $R$. A graded ideal $\mathfrak a\subseteq R$ is \defn{gr-differential} (with respect to $(F,d)_h$) if $d(\mathfrak a)\subseteq\mathfrak a$ and $F(\mathfrak a)\subseteq\mathfrak a$.
\end{definition}

\begin{proposition}\label{O6B9M80}
Let $\phi: R \too S$ be a surjective graded homomorphism with $\ker\phi$ a gr-differential ideal. Then $\phi$ induces
\[
\phi_*: \mathfrak{Der}^{gh}_\Gamma(R) \too \mathfrak{Der}^{gh}_\Gamma(S),
\]
defined by $\phi_*((F_R, d_R)_h) = (F_S, d_S)_h$, where $F_S(\phi(r)) = \phi(F_R(r))$ and $d_S(\phi(r)) = \phi(d_R(r))$.
\end{proposition}

\begin{proof}
Well-definedness follows from the gr-differential property of $\ker(\phi)$: if $\phi(r)=\phi(r')$, then $r-r'\in\ker\phi$, so $F_R(r-r'),d_R(r-r')\in\ker\phi$, and hence $\phi(F_R(r))=\phi(F_R(r'))$, $\phi(d_R(r))=\phi(d_R(r'))$. The derivation properties of $F_S,d_S$ transfer directly from those of $F_R,d_R$ via surjectivity of $\phi$ and the graded homomorphism property.
\end{proof}

Applying Proposition~\ref{O6B9M80} to a graded isomorphism, and separately to its inverse, yields the following.

\begin{corollary}\label{TZ9ATSN}
Let $\varphi: R \too S$ be a graded isomorphism. Then
\[
\Phi_\varphi: \mathfrak{Der}^{gh}_\Gamma(R) \longrightarrow \mathfrak{Der}^{gh}_\Gamma(S),
\quad
\Phi_\varphi((F_R,d_R)_h) = (\varphi \circ F_R \circ \varphi^{-1}, \varphi \circ d_R \circ \varphi^{-1})_h
\]
is a bijection.
\end{corollary}

\begin{proof}
Since $\ker\varphi=0$ is trivially gr-differential, Proposition~\ref{O6B9M80} applies to $\varphi$ and gives $\Phi_\varphi=\varphi_*$; the same reasoning applied to the graded isomorphism $\varphi^{-1}$ gives $\Phi_{\varphi^{-1}}=(\varphi^{-1})_*$. For $(F_R,d_R)_h\in\mathfrak{Der}^{gh}_\Gamma(R)$, the defining formula $F_S(\varphi(r))=\varphi(F_R(r))$ shows that $\Phi_{\varphi^{-1}}(\Phi_\varphi((F_R,d_R)_h))$ agrees with $(F_R,d_R)_h$ on every $r\in R$, and symmetrically for $\Phi_\varphi\circ\Phi_{\varphi^{-1}}$. Hence $\Phi_{\varphi^{-1}}$ is a two-sided inverse of $\Phi_\varphi$, which is therefore a bijection.
\end{proof}

\begin{definition}\label{GHDHOM}
Let $R$ and $S$ be $\Gamma$-graded rings, equipped with generalized homogeneous derivations $(F_R, d_R)_h$ and $(F_S, d_S)_h$, respectively. A graded homomorphism $\phi: R \too S$ is called a \defn{ghd-homomorphism} if it satisfies the compatibility conditions $\phi \circ F_R = F_S \circ \phi$ and $\phi \circ d_R = d_S \circ \phi$, that is, if the following diagrams
\[
\begin{tikzcd}
R \arrow[r, "F_R"] \arrow[d, "\phi"'] & R \arrow[d, "\phi"] \\
S \arrow[r, "F_S"'] & S
\end{tikzcd}
\qquad\text{and}\qquad
\begin{tikzcd}
R \arrow[r, "d_R"] \arrow[d, "\phi"'] & R \arrow[d, "\phi"] \\
S \arrow[r, "d_S"'] & S
\end{tikzcd}
\]
commute.
\end{definition}

\begin{definition}\label{CATGH}
The category $\mathscr{G}_{\Gamma}^{h}$ has, as objects, the triples $(R, F, d)$ where $R$ is $\Gamma$-graded and $(F, d)_h$ is a generalized homogeneous derivation, and, as morphisms, the ghd-homomorphisms between the underlying rings.
\end{definition}

We now turn to the algebraic structure carried by gr-generalized derivations themselves.

\begin{proposition}\label{PMODSTRUCT}
$p\mathfrak{Der}^{gh}_\Gamma(R)$ forms a $Z(R) \cap R_e$-module under pointwise addition and the scalar multiplication $r \cdot (F,d)_h = (rF, rd)_h$.
\end{proposition}

\begin{proof}
$p\mathfrak{Der}^{gh}_\Gamma(R)$ is an additive group under pointwise addition. For $r \in Z(R) \cap R_e$ and $(F,d)_h \in p\mathfrak{Der}^{gh}_\Gamma(R)$, the pair $(rF, rd)_h$ again lies in $p\mathfrak{Der}^{gh}_\Gamma(R)$: the generalized-derivation identity is inherited from the centrality of $r$, and degree preservation follows since $x\in R_\tau$ gives $F(x),d(x)\in R_\tau$, so that $(rF)(x)=r(F(x))\in R_eR_\tau=R_\tau$ and likewise $(rd)(x)\in R_\tau$, using $r\in R_e$.
\end{proof}

\begin{proposition}\label{EPXSK28}
$p\mathfrak{Der}^{gh}_\Gamma(R)$ is a Lie ring\footnote{A Lie ring is an algebraic structure similar to a Lie algebra, but whose underlying additive group is just an abelian group rather than a vector space over a field.} under the bracket
\[
[(F_1,d_1)_h,(F_2,d_2)_h] = (F_1 \circ F_2 - F_2 \circ F_1,\, d_1 \circ d_2 - d_2 \circ d_1)_h,
\]
that is, $p\mathfrak{Der}^{gh}_\Gamma(R)$ is closed under this bracket, which is biadditive, alternating, and satisfies the Jacobi identity.
\end{proposition}
\begin{proof}
Since $d_1,d_2$ are derivations, the mixed terms in $d_1d_2(xy)$ and $d_2d_1(xy)$ cancel upon subtraction, so $[d_1,d_2]$ is again a homogeneous derivation. Expanding $F_1F_2(xy)$ and $F_2F_1(xy)$ via the generalized-derivation identity, the mixed terms cancel likewise, giving $[F_1,F_2](xy)=[F_1,F_2](x)y+x[d_1,d_2](y)$, with $[F_1,F_2]$ homogeneous; hence $([F_1,F_2],[d_1,d_2])_h\in p\mathfrak{Der}^{gh}_\Gamma(R)$, proving closure. Biadditivity and alternation follow immediately from bilinearity of composition, and the Jacobi identity is the standard commutator identity in an associative ring, applied componentwise to $(F_i)$ and $(d_i)$.
\end{proof}

\begin{remark}
The bracket above is not $Z(R) \cap R_e$-bilinear, so $p\mathfrak{Der}^{gh}_\Gamma(R)$ is a Lie ring rather than a Lie algebra over $Z(R) \cap R_e$: its $Z(R) \cap R_e$-module structure and its Lie bracket do not interact compatibly. Take $R = k[t]$ with the trivial grading and $d = \dfrac{\partial}{\partial t}$, and set $F= (\id_R,0)$ and $G = (D,D)$ in $p\mathfrak{Der}^{gh}_\Gamma(R)$. Then $[F,G] = 0$, so bilinearity over $t \in Z(R) \cap R_e$ would force $[tF,G] = t[F,G] = 0$. But for every $f \in R$, $[tF,G](f)= -f$, so $[tF,G] = (-\id_R,0) \neq 0$.
\end{remark}

We turn next to finite products of gr-generalized derivations. Example \ref{PRODFAIL-EX} showed that the coordinatewise construction fails already for generalized homogeneous derivations that merely preserve homogeneity.

\begin{proposition}\label{PRODFIX}
Let $\{R_i\}_{i\in I}$ be a finite collection of $\Gamma_i$-graded rings, graded as in \eqref{PRODGRAD}, and set $R = \prod_{i \in I} R_i$. The assignment
\[
\Psi : \prod_{i \in I} p\mathfrak{Der}^{gh}_{\Gamma_i}(R_i) \too p\mathfrak{Der}^{gh}_{\prod_{i \in I} \Gamma_i}(R),
\quad
\Psi\bigl((F_i,d_i)_h\bigr)_{i \in I} = (F,d)_h,
\]
where $F((r_i)_{i \in I}) = (F_i(r_i))_{i \in I}$ and $d((r_i)_{i \in I}) = (d_i(r_i))_{i \in I}$, is a well-defined injective map.
\end{proposition}

\begin{proof}
Additivity of $F$ and $d$ is immediate from additivity in each coordinate. Let $r=(r_i)_i \in R_{(\tau_i)_i}$, so that $r_i\in (R_i)_{\tau_i}$ for every $i\in I$. Since each $(F_i,d_i)_h$ is gr-generalized, $F_i(r_i)\in (R_i)_{\tau_i}$ and $d_i(r_i)\in(R_i)_{\tau_i}$, hence
\[
F(r) = (F_i(r_i))_i \in \prod_{i\in I}(R_i)_{\tau_i} = R_{(\tau_i)_i}, \quad d(r)\in R_{(\tau_i)_i}
\]
as well. Thus $F$ and $d$ preserve every homogeneous component $R_{(\tau_i)_i}$ of $R$ -- in particular $F(\mathcal H(R))\subseteq\mathcal H(R)$, in contrast with Example \ref{PRODFAIL-EX} -- and the Leibniz rule for $F$, together with the homogeneity of $d$, holds coordinatewise. Hence $(F,d)_h\in p\mathfrak{Der}^{gh}_{\prod_i \Gamma_i}(R)$, and $\Psi$ is well defined.

For injectivity, suppose $\Psi((F_i,d_i)_h)_{i \in I} = \Psi((F_i',d_i')_h)_{i \in I}$. Evaluating both sides at $e_i(r_i) = (0,\ldots,r_i,\ldots,0)$ gives $F_i(r_i) = F_i'(r_i)$ and $d_i(r_i) = d_i'(r_i)$ for every $i \in I$ and every $r_i \in R_i$, so $(F_i,d_i)_h = (F_i',d_i')_h$ for each $i$.
\end{proof}

\begin{remark}\label{NONSURJ}
$\Psi$ is not surjective in general. Let $B$ be any unital ring with nontrivial multiplication, and let $A_0$ be a nonzero ring with zero multiplication. Set $A := B \times A_0$, with $(x,y)(x',y') = (xx',0)$ for $x,x'\in B$, $y,y'\in A_0$, and $R := A \times A$, equipped with the trivial grading. Every element of $R$ is then homogeneous, since a trivial grading has a single homogeneous component; homogeneity preservation and degree preservation therefore coincide for it, so $p\mathfrak{Der}^{gh}_\Gamma(A) = \mathfrak{Der}^{gh}_\Gamma(A)$ here. Define $\pi : A \too A$ by $\pi(x,y) = (0,y)$, and set
\[
F(u,v) := (\pi(v),\,0_A), \quad d = 0,
\]
for $u,v \in A$. For $X=(u,v)$ and $Y=(u',v')$ in $R$, write $u=(x_1,y_1)$, $v=(x_2,y_2)$, $u'=(x_3,y_3)$, $v'=(x_4,y_4)$ with $x_1,x_2,x_3,x_4\in B$ and $y_1,y_2,y_3,y_4\in A_0$. Since $A_0$ has zero multiplication, $vv' = (x_2,y_2)(x_4,y_4) = (x_2x_4,0)$, so $\pi(vv') = 0_A$ and hence $F(XY) = (\pi(vv'),0_A) = 0_R$. Likewise, $\pi(v)u' = (0,y_2)(x_3,y_3) = (0,0) = 0_A$, so $F(X)Y = (\pi(v)u',\,0_A \cdot v') = 0_R$. Thus $F(XY) = F(X)Y$ for all $X,Y \in R$, so $(F,0)_h \in p\mathfrak{Der}^{gh}_\Gamma(R)$. $F$ mixes the two factors, since its output in the first coordinate depends on $v$, the second input coordinate: it is not of the form $(u,v) \mapsto (F_1(u),F_2(v))$ for any pair of maps $F_1,F_2$, so $(F,0)_h$ does not lie in the image of $\Psi$.
\end{remark}

We now package the product construction into the categorical language above.

\begin{definition}\label{PRODCAT}
The category $p\mathscr{G}_\Gamma^{gh}$ has, as objects, the triples $(R,F,d)$ with $R$ a $\Gamma$-graded ring and $(F,d)_h\in p\mathfrak{Der}^{gh}_\Gamma(R)$, and, as morphisms, the ghd-homomorphisms of Definition~\ref{GHDHOM} between the underlying rings.
\end{definition}

To compute finite products inside $p\mathscr{G}_\Gamma^{gh}$, where every factor is graded by the \emph{same} group $\Gamma$, it is more convenient to grade $\prod_{i\in I}R_i$ diagonally by $\Gamma$ rather than by $\Gamma^{|I|}$ as in \eqref{PRODGRAD}.

\begin{lemma}\label{DIAGGRAD}
Let $\{R_i\}_{i\in I}$ be a finite collection of $\Gamma$-graded rings for a fixed group $\Gamma$, and set $R=\prod_{i\in I}R_i$. For $\tau\in \Gamma$, put $R_\tau := \prod_{i\in I}(R_i)_\tau$. Then $R=\bigoplus_{\tau\in \Gamma}R_\tau$ is a $\Gamma$-grading on $R$, and $\Psi$ of Proposition~\ref{PRODFIX}, computed with every $\Gamma_i$ taken equal to $\Gamma$ and composed with the group homomorphism $\Gamma\hookrightarrow \Gamma^{|I|}$, $\tau\mapsto(\tau,\ldots,\tau)$, restricts to a well-defined injective map
\[
\Psi_\Gamma : \prod_{i\in I}p\mathfrak{Der}^{gh}_\Gamma(R_i) \too p\mathfrak{Der}^{gh}_\Gamma(R).
\]
\end{lemma}

\begin{proof}
Given $r=(r_i)_i\in R$, each $r_i=\sum_{\tau\in \Gamma}(r_i)_\tau$ has finite support $S_i\subseteq \Gamma$ with $S=\bigcup_{i\in I}S_i$ finite, and $\rho_\tau:=((r_i)_\tau)_{i\in I}\in R_\tau$ for $\tau\in S$ (where $(r_i)_\tau:=0$ if $\tau\notin S_i$), so that $r=\sum_{\tau\in S}\rho_\tau$. If $\sum_\tau\rho_\tau=0$ with $\rho_\tau\in R_\tau$, projecting to coordinate $i$ gives $\sum_\tau(a_i)_\tau=0$ in $R_i$ for each $i$, so $(a_i)_\tau=0$ for every $\tau,i$ by uniqueness in $R_i$, hence $\rho_\tau=0$ for every $\tau$. Multiplicativity is again coordinatewise. For $(F_i,d_i)_h\in p\mathfrak{Der}^{gh}_\Gamma(R_i)$ and $r=(r_i)_i\in R_\tau$, i.e., $r_i\in(R_i)_\tau$ for every $i$, degree preservation of each $(F_i,d_i)_h$ gives $F_i(r_i),d_i(r_i)\in (R_i)_\tau$, so $F(r):=(F_i(r_i))_i$ and $d(r):=(d_i(r_i))_i$ both lie in $R_\tau$; the Leibniz rule holds coordinatewise as in Proposition~\ref{PRODFIX}. Injectivity is proved exactly as there.
\end{proof}

\begin{proposition}\label{PRODCATFIN}
$p\mathscr{G}_{\Gamma}^{gh}$ admits finite products.
\end{proposition}

\begin{proof}
Let $\{(R_i, F_i, d_i)\}_{i \in I}$ be a finite family of objects in $p\mathscr{G}_{\Gamma}^{gh}$. By Lemma~\ref{DIAGGRAD}, $R = \prod_{i \in I} R_i$, graded diagonally by $\Gamma$, carries a gr-generalized derivation $(F, d)_h \in p\mathfrak{Der}^{gh}_{\Gamma}(R)$ where $F((r_i)_{i \in I}) = (F_i(r_i))_{i \in I}$ and $d((r_i)_{i \in I}) = (d_i(r_i))_{i \in I}$, making $(R, F, d)$ an object of $p\mathscr{G}_{\Gamma}^{gh}$. To verify the universal property of products, let $(T, F_T, d_T)$ be an arbitrary object of $p\mathscr{G}_{\Gamma}^{gh}$ and let $\{\phi_i: T \too R_i\}_{i \in I}$ be a family of morphisms in $p\mathscr{G}_{\Gamma}^{gh}$. Define $\phi: T \too R$ by $\phi(t) = (\phi_i(t))_{i \in I}$; then $\phi$ is graded, since each $\phi_i$ is, and $\pi_i \circ \phi = \phi_i$ for each $i$, where $\pi_i: R \too R_i$ is the canonical projection, itself a graded ring homomorphism for the diagonal grading. To see that $\phi$ is a morphism in $p\mathscr{G}_{\Gamma}^{gh}$, we check compatibility with the gr-generalized derivations: for any $t \in T$,
\begin{align*}
    \phi(F_T(t)) &= (\phi_i(F_T(t)))_{i \in I}
    = (F_i(\phi_i(t)))_{i \in I}
    = F((\phi_i(t))_{i \in I})
    = F(\phi(t)),
\end{align*}
using that each $\phi_i$ is a ghd-homomorphism, and similarly
\[
\phi(d_T(t)) = (\phi_i(d_T(t)))_{i \in I} = (d_i(\phi_i(t)))_{i \in I} = d((\phi_i(t))_{i \in I}) = d(\phi(t)).
\]
If $\psi: T \too R$ is another morphism in $p\mathscr{G}_{\Gamma}^{gh}$ with $\pi_i \circ \psi = \phi_i$ for each $i$, then for any $t \in T$,
\[
\psi(t) = (\pi_i(\psi(t)))_{i \in I} = (\phi_i(t))_{i \in I} = \phi(t),
\]
so $\psi = \phi$, proving uniqueness.
\end{proof}

We close this section with two further constructions: extension of a gr-generalized derivation to a polynomial variable, and its behavior under tensor product.

\begin{proposition}\label{T7A9TGI}
For a $\Gamma$-graded $R$ with $R[t]$ graded by $\deg(t) = e$, there exists a natural injection
\[
p\mathfrak{Der}^{gh}_\Gamma(R) \hookrightarrow p\mathfrak{Der}^{gh}_\Gamma(R[t])
\]
given by $(F,d)_h \mapsto (F',d')_h$, where $F'(\sum r_i t^i) = \sum F(r_i)t^i$ and $d'(\sum r_i t^i) = \sum d(r_i)t^i$.
\end{proposition}

\begin{proof}
Let $(F,d)_h \in p\mathfrak{Der}^{gh}_\Gamma(R)$ and define $F',d' : R[t] \too R[t]$ by
\[
F'\!\left(\sum_{i=0}^n r_i t^i\right) = \sum_{i=0}^n F(r_i)t^i,
\quad
d'\!\left(\sum_{i=0}^n r_i t^i\right) = \sum_{i=0}^n d(r_i)t^i.
\]
That $(F',d')_h \in p\mathfrak{Der}^{gh}_\Gamma(R[t])$ follows directly, since $F$ and $d$ act termwise. For homogeneity, if $f(t) = \sum_{i=0}^n r_i t^i \in R[t]_\tau$, then $r_i \in R_\tau$, so $F(r_i), d(r_i) \in R_\tau$, and $F'(f), d'(f) \in R[t]_\tau$. For injectivity, if $(F,d)_h \neq (0,0)_h$, then $F(r)\neq 0$ or $d(r)\neq0$ for some $r\in R$, so $F'(r)\ne0$ or $d'(r)\ne0$ when $r$ is viewed as a constant polynomial.
\end{proof}

The injection in Proposition~\ref{T7A9TGI} is not surjective in general. For instance, take $R = \mathbb{C}$ with the trivial grading and endow $R[t]$ with the standard $\mathbb{Z}$-grading. Define a derivation $d : R[t] \too R[t]$ by
\[
d\!\left(\sum_{i} a_i t^i\right) = \sum_{i} i\,a_i t^{i},
\]
so that $d(t)=t \neq 0$ and $d$ is homogeneous of degree $0$. Then $(d,d)_h \in p\mathfrak{Der}^{gh}_{\mathbb{Z}}(R[t])$, but $(d,d)_h$ cannot belong to the image of the injection, since any element of the image satisfies $d'(t)=0$.

\begin{proposition}\label{4KFS6VF}
For $\Gamma$-graded $k$-algebras $R, S$ with $(F_R,d_R)_h \in p\mathfrak{Der}^{gh}_\Gamma(R)$ and $(F_S,d_S)_h \in p\mathfrak{Der}^{gh}_\Gamma(S)$, define
\[
F_{R \otimes S}(r \otimes s) = F_R(r) \otimes s + r \otimes F_S(s), \quad
d_{R \otimes S}(r \otimes s) = d_R(r) \otimes s + r \otimes d_S(s).
\]
Then $(F_{R \otimes S},d_{R \otimes S})_h \in p\mathfrak{Der}^{gh}_\Gamma(R \otimes_k S)$.
\end{proposition}

\begin{proof}
Extend $F_{R\otimes S}$ and $d_{R\otimes S}$ to all of $R\otimes_k S$ by linearity. By linearity it suffices to check the generalized-derivation identity on homogeneous tensors $u=r_1\otimes s_1$, $v=r_2\otimes s_2$:
\begin{align*}
F_{R \otimes S}(uv) &= F_{R \otimes S}(r_1r_2 \otimes s_1s_2)
= F_R(r_1r_2) \otimes s_1s_2 + r_1r_2 \otimes F_S(s_1s_2) \\
&= (F_R(r_1)r_2 + r_1d_R(r_2)) \otimes s_1s_2 + r_1r_2 \otimes (F_S(s_1)s_2 + s_1d_S(s_2)) \\
&= \bigl[F_R(r_1) \otimes s_1 + r_1 \otimes F_S(s_1)\bigr](r_2 \otimes s_2)+(r_1 \otimes s_1)\bigl[d_R(r_2) \otimes s_2+r_2 \otimes d_S(s_2)\bigr] \\
&= F_{R \otimes S}(u)\,v + u\,d_{R \otimes S}(v).
\end{align*}
The same computation, with $F_R,F_S$ replaced by $d_R,d_S$, gives the Leibniz rule for $d_{R\otimes S}$:
\[
d_{R\otimes S}(uv) = d_{R\otimes S}(u)v + u\,d_{R\otimes S}(v).
\]
For homogeneity, if $r\in R_\tau$ and $s\in S_\sigma$, then $r\otimes s\in(R\otimes_K S)_{\tau\sigma}$, and since $F_R(r)\in R_\tau$, $F_S(s)\in S_\sigma$, $d_R(r)\in R_\tau$, $d_S(s)\in S_\sigma$, both $F_{R\otimes S}(r\otimes s)$ and $d_{R\otimes S}(r\otimes s)$ lie in $(R\otimes_K S)_{\tau\sigma}$. Hence $(F_{R\otimes S},d_{R\otimes S})_h\in p\mathfrak{Der}^{gh}_\Gamma(R\otimes_K S)$.
\end{proof}

\section{Generalized homogeneous derivations on gr-prime rings}\label{9J1ICTX}

We now specialize to gr-prime rings. Two structural facts are recorded first, before we turn to the commutativity criteria that are the main results of the paper.

\begin{proposition}\label{PROP-FNONZERO}
Let $R$ be a gr-prime ring and $(F,d)_h$ a generalized homogeneous derivation of $R$. If $d \neq 0$, then $F \neq 0$.
\end{proposition}

\begin{proof}
Assume $F = 0$. The relation $F(xy) = F(x)y + xd(y)$ then yields $xd(y) = 0$ for all $x,y \in R$, so $xRd(y) = 0$ for all $x,y\in R$. For any nonzero $r \in \mathcal{H}(R)$, $rRd(y) = 0$ for all $y \in R$, so by \cite[Proposition 2.1]{yassine}, $d(y) = 0$ for all $y \in R$, contradicting $d \neq 0$.
\end{proof}

\begin{theorem}\label{THM-DECOMP}
For gr-prime rings $R$, there exists a canonical decomposition
\[
p\mathfrak{Der}^{gh}_\Gamma(R) = p\mathfrak{Der}^{h}_\Gamma(R) \oplus \mathcal{C}_\Gamma(R),
\]
where $\mathcal{C}_\Gamma(R) = \{ F \in p\mathfrak{Der}^{gh}_\Gamma(R) \, :\, F \text{ has zero associated derivation} \}$.
\end{theorem}

\begin{proof}
For $(F,d)_h \in p\mathfrak{Der}^{gh}_\Gamma(R)$, write $F = d + (F-d)$ and set $F_1 = d$, $F_2 = F-d$. Then $F_1 \in p\mathfrak{Der}^{h}_\Gamma(R)$, and $F_2(xy) = F_2(x)y$ for all $x,y \in R$, with degree preservation inherited from $F$ and $d$, so $F_2 \in \mathcal{C}_\Gamma(R)$. It remains to show $p\mathfrak{Der}^{h}_\Gamma(R) \cap \mathcal{C}_\Gamma(R) = 0$.

Let $H$ lie in this intersection. The derivation property gives $H(xy) = H(x)y + xH(y)$, while $H \in \mathcal{C}_\Gamma(R)$ gives $H(xy) = H(x)y$. Hence $xH(y) = 0$ for all $x,y \in R$, so $xRH(y) = 0$ for all $x,y \in R$. For any $r \in \mathcal{H}(R)\setminus\{0\}$, $rRH(y) = 0$ for all $y$, and \cite[Proposition 2.1]{yassine} gives $H(y) = 0$ for all $y \in R$, so $H \equiv 0$.
\end{proof}

\subsection{Some commutativity criteria}\label{G1WU16H}

We begin with results for plain homogeneous derivations, before turning to their generalized counterparts.

\begin{proposition}\label{HPCF07U}
Let $R$ be a gr-prime ring and $\mathfrak{a}$ a nonzero graded ideal of $R$ such that
\[
[x,y]\in Z(R)
\quad\text{or}\quad
x\circ y\in Z(R)
\]
for all $x,y\in \mathfrak{a}$. Then $R$ is a commutative graded ring.
\end{proposition}

The proof uses the following lemma.

\begin{lemma}\label{MSND8DU}
Let $R$ be a gr-prime ring. Then:
\begin{enumerate}
\item[(1)] if $\mathfrak{a}$ is a nonzero graded ideal of $R$ and $a\mathfrak{a}b=0$ where $a \in \mathcal{H}(R)$ or $b\in \mathcal{H}(R)$, then $a=0$ or $b=0$;
\item[(2)] if $d$ is a homogeneous derivation of $R$ and $ad(x) = 0$ or $d(x)a = 0$ for all $x \in R$, then $a = 0$ or $d = 0$.
\end{enumerate}
\end{lemma}

\begin{proof}
(1) Let $a=\sum_{\tau\in \Gamma}a_\tau \in R$ and $b \in \mathcal{H}(R)\setminus\{0\}$ with $a\mathfrak{a}b=0$. For every $r \in \mathfrak{a}\cap \mathcal{H}(R)$, $arb=0$, so $\sum_{\tau\in \Gamma}a_\tau rb=0$; since the terms $a_\tau rb$ lie in pairwise distinct homogeneous components (the map $\tau\mapsto \tau\deg(r)\deg(b)$ is injective), uniqueness of the homogeneous decomposition gives $a_\tau rb=0$ for every $\tau \in \Gamma$ and every $r \in \mathfrak{a} \cap \mathcal{H}(R)$. As $\mathfrak a$ is an ideal, $rx\in\mathfrak a\cap\mathcal H(R)$ for every homogeneous $x\in R$, so $a_\tau (rx) b=0$ as well; by additivity in $x$, $a_\tau r x b=0$ for every $x\in R$, hence $a_\tau \mathfrak{a}Rb=0$ for all $\tau$. By \cite[Proposition 2.1]{yassine}, $a_\tau \mathfrak{a}=0$, hence $a_\tau R\mathfrak{a}=0$, for all $\tau$, so $a_\tau=0$ for all $\tau$ and $a=0$.

(2) Suppose $ad(x) = 0$ and $a\neq 0$. Replacing $x$ by $xy$ gives $axd(y)=0$ for all $x,y\in R$, so $aRd(x)=0$ for all $x\in R$, hence $aRd(r)=0$ for all $r\in \mathcal{H}(R)$. By \cite[Proposition 2.1]{yassine}, $d(r)=0$ for all $r\in \mathcal{H}(R)$, so $d=0$.
\end{proof}

\begin{proof}[Proof of Proposition \ref{HPCF07U}]
Suppose first that $[x,y] \in Z(R)$ for all $x,y\in \mathfrak{a}$. Then
\begin{equation}\label{O2RAK1V}
[z,[x,y]]=0
\end{equation}
for all $x,y\in \mathfrak{a}$, $z\in R$. Replacing $y$ by $yx$ in \eqref{O2RAK1V} and simplifying gives
\begin{equation}\label{GWWQC38}
[x,y][z,x]=0
\end{equation}
for all $x,y\in \mathfrak{a}$, $z\in R$. Substituting $zy$ for $z$ in \eqref{GWWQC38} yields $[x,y]z[x,y]=0$, so $[x,y]R[x,y]=0$ for all $x,y\in \mathfrak{a}$; since $\mathfrak{a}$ is graded, $[r_1,r_2]R[r_1,r_2]=0$ for all $r_1, r_2\in \mathfrak{a}\cap \mathcal{H}(R)$. Gr-primeness gives $[r_1,r_2]=0$, hence $[x,y]=0$ for all $x,y\in \mathfrak{a}$: $\mathfrak{a}$ is commutative, and $R$ is commutative by \cite[Proposition 2.1]{yassine}.

Suppose instead that $x\circ y\in Z(R)$ for all $x,y\in \mathfrak{a}$, so $[x\circ y,z]=0$ for all $x,y\in \mathfrak{a}$, $z\in R$. Replacing $y$ by $yx$ and simplifying gives
\begin{equation}\label{HYF9SXM}
(x\circ y)[x,z]=0
\end{equation}
for all $x,y\in \mathfrak{a}$, $z\in R$. Substituting $sz$ for $z$ gives $(r_1\circ r_2)R[r_1,z]=0$ for all $r_1,r_2\in \mathfrak{a}\cap \mathcal{H}(R)$, $z\in R$, so by \cite[Proposition 2.1]{yassine}, $x\circ y=0$ or $[x,z]=0$ for all $x,y\in \mathfrak{a}$, $z\in R$. If $[x,z]=0$ throughout, $\mathfrak{a}$ is a central graded ideal and $R$ is commutative. Otherwise $x\circ y=0$ for all $x,y\in\mathfrak a$; replacing $y$ by $yz$ gives $y[x,z]=0$, and Lemma  \ref{MSND8DU}(2) again gives $[x,z]=0$ for all $x\in \mathfrak{a}$, $z\in R$. In both cases $\mathfrak{a}$ is a central graded ideal, so $R$ is commutative.
\end{proof}

The next result characterizes when compositions of homogeneous derivations force commutativity.

\begin{theorem}\label{HT6T0UZ}
Let $R$ be a gr-prime ring of characteristic different from $2$. Suppose $d_{1}$ and $d_{2}$ are nonzero homogeneous derivations of $R$ such that
\[
d_{1}d_{2}(x) \in Z(R)
\]
for all $x\in R$. Then $R$ is a commutative graded ring.
\end{theorem}

\begin{proof}
By hypothesis,
\begin{equation}\label{6PR5N8I}
d_{1}d_{2}(x)\in Z(R)
\end{equation}
for all $x\in R$. Replacing $x$ by $[x,y]$ in \eqref{6PR5N8I} and expanding, using that $d_1,d_2$ are derivations so $d_i([x,y])=[d_i(x),y]+[x,d_i(y)]$, and that $[d_1d_2(x),y]=[x,d_1d_2(y)]=0$ since $d_1d_2(x),d_1d_2(y)\in Z(R)$ by \eqref{6PR5N8I}, gives
\begin{equation}\label{4P40MGG}
[d_2(x),d_{1}(y)]+[d_{1}(x),d_{2}(y)]\in Z(R)
\end{equation}
for all $x,y\in R$. Putting $y=d_{2}(z)$ in \eqref{4P40MGG}, and using again that $[d_2(x),d_1d_2(z)]=0$ since $d_1d_2(z)\in Z(R)$, yields $[d_1(x),d_{2}^{2}(z)]\in Z(R)$ for all $x,z\in R$, so in particular $[d_{2}^{2}(r),d_{1}(y)]\in Z(R)$ for all $r\in \mathcal{H}(R)$, $y\in R$. By \cite[Lemma 2.2]{yassine}, either $d_{2}^{2}(r)\in Z(R)$ for all $r\in \mathcal{H}(R)$, or $d_{1}=0$; the latter is excluded by hypothesis, so $d_{2}^{2}(x)\in Z(R)$ for all $x\in R$ by additivity. Taking $[x,z]$ instead of $x$ gives $2[d_{2}(x),d_{2}(z)]\in Z(R)$ for all $x,z\in R$, and since $\operatorname{char}R \neq 2$, $[d_{2}(x),d_{2}(z)]\in Z(R)$ for all $x,z\in R$. By \cite[Theorem 3.5]{yassine}, $R$ is commutative.
\end{proof}

The next example shows that the gr-primeness hypothesis in Theorem~\ref{HT6T0UZ} cannot be weakened to gr-semiprimeness.

\begin{example}\label{EX-HT6T0UZ}
Consider $R=\mathbb{C}[t_1,t_2,t_3,t_4]\times M_{2}(\mathbb{C})$, where $\mathbb{C}[t_1,t_2,t_3,t_4]$ carries the standard $\mathbb{Z}$-grading and $M_2(\mathbb{C})$ the elementary $\mathbb{Z}_4$-grading with $(M_2(\mathbb{C}))_0=\{\mathrm{diag}(a,b)\}$, $(M_2(\mathbb{C}))_2=\{ae_{12}+be_{21}\}$, $(M_2(\mathbb{C}))_1=(M_2(\mathbb{C}))_3=0$, and $R$ carries the product $\mathbb{Z}\times \mathbb{Z}_4$-grading. Since $(1,0)$ and $(0,I_2)$ are nonzero homogeneous elements of degree $(0,0)$ with $(1,0)R(0,I_2)=0$, $R$ is not gr-prime; it is gr-semiprime. Define homogeneous derivations $d_1, d_2: R \too R$ by
\[
d_1\left(f,M\right)=\left(t_2t_4\frac{\partial f}{\partial t_1},0\right),
\quad
d_2\left(f, M\right)=\left(t_1t_3\frac{\partial f}{\partial t_2},0\right).
\]
Then $d_1d_2(x)\in Z(R)$ for all $x\in R$, so the hypothesis of Theorem \ref{HT6T0UZ} holds, yet $R$ is noncommutative (as $M_2(\mathbb{C})$ is a noncommutative factor).
\end{example}

In \cite{achraf2}, it was shown that a prime ring $R$ with a nonzero ideal $\mathfrak{a}$ is commutative if it admits a generalized derivation $F$ satisfying $F(xy)\pm xy\in Z(R)$ or $F(x)F(y)\pm xy\in Z(R)$ for all $x,y\in \mathfrak{a}$. We extend this result to gr-prime rings for generalized homogeneous derivations. The graded setting forces a genuinely new step: gr-primeness only sees homogeneous elements, so the case distinctions that close the argument over a prime ring must first be established on $\mathfrak{a}\cap\mathcal H(R)$, and then shown not to split differently across different homogeneous degrees. We isolate this step as a single lemma, used three times below.

\begin{lemma}\label{GRDICH}
Let $R$ be a gr-prime ring, $\mathfrak a\ne\{0\}$ a graded ideal of $R$, and $d$ a homogeneous derivation of $R$. Let $S=\bigoplus_\tau(S\cap R_\tau)$ be a graded additive subgroup of $R$, and let $(c_x)_{x\in E}$, $E\ne\emptyset$, be additive maps $c_x:S\to R$ satisfying
\[
c_x(z)\,R\,\mathfrak a\,d(z)=0\quad\text{for all }z\in S,\ x\in E.
\]
Set $Z_1=\{z\in S:c_x(z)=0\ \forall x\in E\}$ and $Z_2=\{z\in S:d(z)=0\}$. Then $S=Z_1$ or $S=Z_2$.
\end{lemma}

\begin{proof}
For $w\in\mathcal H(R)$ we use the following two facts.

(i) If $\mathfrak a\,d(w)=0$ then $d(w)=0$: pick $0\ne r\in\mathfrak a\cap\mathcal H(R)$; then $rRd(w)=0$, and gr-primeness gives $d(w)=0$.

(ii) If $\mathfrak a\,d(w)\ne0$, some homogeneous $r \in\mathfrak a$ satisfies $rd(w)\ne0$: write $a=\sum_\tau a_\tau\in\mathfrak a$ with $ad(w)\ne0$; the components $a_\tau d(w)$ lie in distinct graded components (since the map $\tau\longmapsto \tau \deg(d(w))$ is injective), so $a_{\sigma}d(w)\ne0$ for some $\sigma$, and we set $r=a_{\sigma}$.

Let $z\in S\cap\mathcal H(R)$ with $d(z)\ne0$. By (i), $\mathfrak a\,d(z)\ne0$, so (ii) gives $r \in\mathfrak a\cap \mathcal{H}(R)$ with $r d(z)\ne0$. The hypothesis then gives $c_x(z)R(r d(z))=0$ for every $x$, and gr-primeness (\cite[Proposition 2.1]{yassine}) forces $c_x(z)=0$ for every $x$. Thus $z\in Z_1$, and $S\cap\mathcal H(R)\subseteq Z_1\cup Z_2$.

Suppose $S\cap\mathcal H(R)$ is contained in neither $Z_1$ nor $Z_2$. Then there is $u\in(S\cap\mathcal H(R))\setminus Z_1$, and since $S\cap\mathcal H(R)\subseteq Z_1\cup Z_2$, this forces $u\in Z_2\setminus Z_1$, so $d(u)=0$ and $c_{x_0}(u)\ne0$ for some $x_0$. Likewise there is $v\in(S\cap\mathcal H(R))\setminus Z_2\subseteq Z_1\setminus Z_2$, so $c_x(v)=0$ for all $x$, and $d(v)\ne0$.

Applying the hypothesis to $z=u+v\in S$ and using $d(u)=0$, $c_x(v)=0$ reduces $c_x(u+v)R\mathfrak a\,d(u+v)=0$ to
\[
c_x(u)\,R\,\mathfrak a\,d(v)=0\quad\text{for all }x.
\]
Since $d(v)\ne0$, (i) gives $\mathfrak a\,d(v)\ne0$, so (ii) gives $r\in\mathfrak a\cap\mathcal H(R)$ with $rd(v)\ne0$; then $c_x(u)R(rd(v))=0$ for all $x$, and gr-primeness forces $c_x(u)=0$ for all $x$, in particular $c_{x_0}(u)=0$, contradicting $c_{x_0}(u)\ne0$. Hence $S\cap\mathcal H(R)$ is contained in $Z_1$ or in $Z_2$.

Say the former holds. Since $S=\bigoplus_\tau(S\cap R_\tau)$, every $z\in S$ is a finite sum of its homogeneous components $z_\tau\in S\cap R_\tau\subseteq Z_1$, and $Z_1$ is additive, so $z\in Z_1$, hence $S=Z_1$. The other case gives $S=Z_2$ the same way.
\end{proof}

\begin{theorem}\label{S9IM0YL}
Let $R$ be a gr-prime ring and $\mathfrak{a}$ a nonzero graded ideal of $R$. If $R$ admits a generalized homogeneous derivation $F$ with associated nonzero homogeneous derivation $d$ such that
\[
F(xy)\pm xy\in Z(R)
\]
for all $x,y\in \mathfrak{a}$, then $R$ is commutative.
\end{theorem}

\begin{proof}
Consider $F(xy)-xy\in Z(R)$ for all $x,y\in \mathfrak{a}$. As in the proof of \cite[Theorem~2.1]{achraf2}, this yields $[z,z_{1}]xyd(z)=0$ for all $x,y,z, z_{1}\in \mathfrak{a}$, hence
\[
[z,z_1]x\,R\,\mathfrak{a}\,d(z) = 0 \quad \text{for all } x,z_1,z\in\mathfrak a,
\]
an identity valid on all of $\mathfrak a$. Apply Lemma \ref{GRDICH} with $S=\mathfrak a$, index set $E = \mathfrak a\times\mathfrak a$, and $c_{(x,z_1)}(z) = [z,z_1]x$, additive in $z$ since the commutator is biadditive. This gives $\mathfrak a = Z_1$ or $\mathfrak a = Z_2$, where $Z_1 = \{z\in\mathfrak a : [z,z_1]x=0\ \forall x,z_1\in\mathfrak a\}$ and $Z_2 = \{z\in\mathfrak a : d(z) = 0\}$.

If $\mathfrak a = Z_1$, then $[z,z_1]\mathfrak a = 0$ for all $z,z_1\in\mathfrak a$; for $r\in\mathfrak a$, $rs\in\mathfrak a$ for every $s\in R$ (as $\mathfrak a$ is an ideal), so $[z,z_1](rs) = 0$, giving $[z,z_1]Rr = 0$ for every $r\in\mathfrak a$, in particular $[z,z_1]Rr_0=0$ for some nonzero $r_0\in\mathfrak a\cap\mathcal H(R)$. Lemma   \ref{MSND8DU}(1) gives $[z,z_1]=0$ for all $z,z_1\in\mathfrak a$: $\mathfrak a$ is commutative, and $R$ is commutative by \cite[Proposition 2.1]{yassine}.

If instead $\mathfrak a = Z_2$, then $d(z)=0$ for all $z\in\mathfrak a$, so $d=0$ on $R$ by \cite[Lemma 2.6]{yassine}, contradicting $d\neq 0$.

The case $F(xy)+xy\in Z(R)$ reduces to the first by considering $-F$ instead of $F$.
\end{proof}

\begin{theorem}\label{YMQC1A7}
Let $R$ be a gr-prime ring and $\mathfrak{a}$ a nonzero graded ideal of $R$. If $R$ admits two generalized homogeneous derivations $F_1$ and $F_2$ with associated nonzero homogeneous derivations $d_1$ and $d_2$, respectively, such that
\[
F_1(x)F_2(y)\pm xy\in Z(R)
\]
for all $x,y\in \mathfrak{a}$, then $R$ is commutative.
\end{theorem}

\begin{proof}
Consider
\begin{equation}\label{RVHS1K3}
F_1(x)F_2(y)-xy\in Z(R)
\end{equation}
for all $x,y \in \mathfrak{a}$. Substituting $yz$ for $y$ in \eqref{RVHS1K3} gives
\begin{equation}\label{AQ6L6GT}
\bigl(F_1(x)F_2(y)-xy\bigr)z+F_1(x)yd_2(z)\in Z(R)
\end{equation}
for all $x,y\in \mathfrak{a}$, $z\in R$. Taking the commutator of \eqref{AQ6L6GT} with $z$ gives
\begin{equation}\label{6MJXK87}
F_1(x)[y d_2(z),z]+[F_1(x),z]y d_2(z)=0,
\end{equation}
and substituting $F_1(x)y$ for $y$ in \eqref{6MJXK87} gives
\[
F_1(x)[F_1(x)yd_2(z),z]+[F_1(x),z]F_1(x)yd_2(z)=0.
\]
Expanding the inner commutator by $[F_1(x)yd_2(z),z]=F_1(x)[yd_2(z),z]+[F_1(x),z]yd_2(z)$ turns this into
\[
F_1(x)^2[yd_2(z),z]+F_1(x)[F_1(x),z]yd_2(z)+[F_1(x),z]F_1(x)yd_2(z)=0,
\]
and subtracting $F_1(x)$ times \eqref{6MJXK87} itself, namely $F_1(x)^2[yd_2(z),z]+F_1(x)[F_1(x),z]yd_2(z)=0$, cancels the first two terms, leaving
\begin{equation}\label{OL73D9P}
[F_1(x),z]F_1(x) y d_2(z)=0
\end{equation}
for all $x,y\in \mathfrak{a}$, $z\in R$, an identity valid for every $z\in R$, homogeneous or not. Since $y$ ranges over the ideal $\mathfrak a$, \eqref{OL73D9P} yields
\[
[F_1(x),z]F_1(x)\,R\,\mathfrak{a}\,d_2(z) = 0\quad\text{for all } x\in\mathfrak a,\ z\in R.
\]
Apply Lemma \ref{GRDICH} with $S = R$, $d = d_2$, index set $E=\mathfrak a$, and $c_x(z) = [F_1(x),z]F_1(x)$, additive in $z$ since the commutator is additive in its second argument. This gives $R = Z_1$ or $R = Z_2$, where $Z_1=\{z\in R : [F_1(x),z]F_1(x)=0\ \forall x\in\mathfrak a\}$ and $Z_2=\{z\in R : d_2(z)=0\}$.

If $R=Z_2$, then $d_2\equiv 0$, contradicting $d_2\neq 0$. Hence $R=Z_1$: $[F_1(x),z]F_1(x)=0$ for all $x\in\mathfrak a$, $z\in R$. Replacing $z$ by $zz'$ gives $[F_1(x),z]RF_1(x)=0$, so by \cite[Proposition 2.1]{yassine}, $F_1(r)=0$ or $[F_1(r),z]=0$ for all $z\in R$; in both cases $[F_1(x),z]=0$ for all $x\in\mathfrak a$, $z\in R$. Replacing $x$ by $xz$ gives
\begin{equation}\label{AU5BHXT}
x[d_1(z),z]+[x,z]d_1(z)=0
\end{equation}
for all $x\in \mathfrak{a}$, $z\in R$, again for every $z\in R$. Substituting $sx$ for $x$ and using \eqref{AU5BHXT} to simplify $[sx,z]=s[x,z]+[s,z]x$ gives
\[
[s,z]xd_1(z)=0 \quad \text{for all } s\in R,\ x\in\mathfrak a,\ z\in R,
\]
and, as $x$ ranges over the ideal $\mathfrak a$, this gives $[s,z]R\mathfrak a\,d_1(z) = 0$ for all $s,z\in R$.

Apply Lemma \ref{GRDICH} once more, with $S=R$, $d=d_1$, index set $E=R$, and $c_s(z)=[s,z]$, additive in $z$. This gives $R=Z_1'$ or $R=Z_2'$, where $Z_1'=\{z\in R:[s,z]=0\ \forall s\in R\}=Z(R)$ and $Z_2'=\{z\in R:d_1(z)=0\}$. If $R=Z_2'$, then $d_1\equiv0$, contradicting $d_1\neq0$; hence $R=Z_1'=Z(R)$, and $R$ is commutative.

The case $F_1(x)F_2(y)+xy\in Z(R)$ reduces to the above by considering $-F_1$ instead of $F_1$.
\end{proof}

The following example shows that the gr-primeness hypothesis in Theorems \ref{S9IM0YL} and \ref{YMQC1A7} cannot be omitted.

\begin{example}\label{EX-FINAL}
Let
\[
R=\mathbb{C}[t_1,t_2,t_3]\times \left\{\begin{pmatrix} a & b \\ 0& 0 \end{pmatrix}\;:\; a,b\in \mathbb{C} \right\}
\]
with the $\mathbb{Z}\times \mathbb{Z}_2$-grading in which $\mathbb{C}[t_1,t_2,t_3]$ carries the standard $\mathbb{Z}$-grading and the second factor is graded by $\deg\begin{pmatrix}a&0\\0&0\end{pmatrix}=0$, $\deg\begin{pmatrix}0&b\\0&0\end{pmatrix}=1$; then $R$ is not gr-prime, since $(1,0)$ and $\left(0,\begin{smallmatrix}0&1\\0&0\end{smallmatrix}\right)$ are nonzero homogeneous elements with product $0$. Let $\mathfrak{a}=\mathbb{C}[t_1,t_2,t_3]\times \left\{ \begin{pmatrix} 0 & a \\ 0 & 0 \end{pmatrix} \;:\; a \in \mathbb{C} \right\}$, a nonzero graded ideal of $R$, and define
\[
F_1(f, M) = \Bigl(t_3\Bigl(f+\frac{\partial f}{\partial t_3}\Bigr),0\Bigr),
\quad
F_2(f, M) = d_2(f,M) = \Bigl(t_1\frac{\partial f}{\partial t_2},0\Bigr),
\quad
d_1(f, M) = \Bigl(t_2t_3\frac{\partial f}{\partial t_3},0\Bigr).
\]
Then $(F_1, d_1)_h$ and $(F_2, d_2)_h$ are generalized homogeneous derivations on $R$ satisfying $F_1(xy) \pm xy \in Z(R)$ and $F_1(x)F_2(y) \pm xy \in Z(R)$ for all $x, y \in \mathfrak{a}$, yet $R$ is noncommutative.
\end{example}

\section{Generalized homogeneous derivations on graded modules}\label{YC4OM3W}

We now extend the theory to graded modules, introducing generalized homogeneous derivations on modules, their functorial behavior, and the associated category.

\begin{definition}\label{IP98Q1Z}
Let $R$ be a $\Gamma$-graded ring and $M$ a $\Gamma$-graded $R$-module. An additive mapping $F_M : M \too M$ is a \defn{generalized homogeneous derivation} if there exists a homogeneous derivation $d : R \too R$ such that:
\begin{enumerate}
\item[(1)] $F_M(rm) = d(r)m + rF_M(m)$ for all $r \in R$, $m \in M$;
\item[(2)] $F_M(m) \in \mathcal{H}(M)$ for all $m \in \mathcal{H}(M)$.
\end{enumerate}
We denote such pairs by $(F_M,d)_{h,M}$ and write $\mathfrak{Der}^{gh}_\Gamma(R,M)$ for the set of all generalized homogeneous derivations on $M$.
\end{definition}

\begin{example}\label{X272HPO}
Let $R = \mathbb{C}[t_1,t_2]$ with the standard $\mathbb{Z}$-grading, and let $M = R^2$ with grading $M_n = \{(f_1,f_2) \,:\, f_i \in R_n\}$. Define
\[
F(f_1,f_2) = \left(\frac{\partial f_1}{\partial t_1}, \frac{\partial f_2}{\partial t_1}\right)
\]
with associated derivation $d(f) = \frac{\partial f}{\partial t_1}$. Then $(F,d)_{h,M} \in \mathfrak{Der}^{gh}_\Gamma(R,M)$.
\end{example}

\begin{definition}\label{XOI65FD}
A graded submodule $N \subseteq M$ is \defn{gr-differential} with respect to $(F_M,d)_{h,M}$ if $F_M(N) \subseteq N$.
\end{definition}

\begin{example}\label{EX-XOI65FD}
The generalized homogeneous derivation $(F_M,d)_{h,M}$ of Example~\ref{X272HPO} makes the graded submodule $N = 0 \oplus R \subseteq M$ gr-differential.
\end{example}

\begin{definition}\label{UN5Z3JQ}
A generalized homogeneous derivation $(F_M,d)_{h,M}$ is \defn{gr-generalized} if $F_M(M_\tau) \subseteq M_\tau$ and $d(R_\tau) \subseteq R_\tau$ for all $\tau \in \Gamma$. The set of gr-generalized derivations on $M$ is denoted $p\mathfrak{Der}^{gh}_\Gamma(R,M)$.
\end{definition}

\begin{proposition}\label{1U322N6}
$p\mathfrak{Der}^{gh}_\Gamma(R,M)$ forms a $Z(R) \cap R_e$-module under pointwise operations and the scalar multiplication $a \cdot (F_M,d)_{h,M} = (aF_M, ad)_{h,M}$ for $a \in Z(R) \cap R_e$.
\end{proposition}

\begin{proof}
Centrality of scalars gives
\[
(aF_M)(rm) = aF_M(rm) = a(d(r)m + rF_M(m)) = (ad)(r)m + r(aF_M)(m),
\]
and degree preservation follows from $a \in R_e$ together with the grading properties of $F_M$ and $d$.
\end{proof}

\begin{proposition}\label{E8N6QIW}
For finite families $\{M_i\}_{i \in I}$ of graded $R$-modules:
\begin{enumerate}
\item if $(F_{M_i},d)_{h,M_i} \in p\mathfrak{Der}^{gh}_\Gamma(R,M_i)$ share the same associated derivation $d$, then
\[
F_{\bigoplus M_i}((m_i)_i) = (F_{M_i}(m_i))_i
\]
defines a canonical gr-generalized derivation on $\bigoplus_{i \in I} M_i$;
\item if, in addition, $R$ is commutative and $(F_M,d)_{h,M}$, $(F_N,d)_{h,N}$ have the same associated derivation $d$, then
\[
F_{M \otimes N}(m \otimes n) = F_M(m) \otimes n + m \otimes F_N(n)
\]
defines a canonical gr-generalized derivation on $M \otimes_R N$.
\end{enumerate}
\end{proposition}

\begin{proof}
(1) For $(m_i)_{i \in I} \in \bigoplus_{i \in I} M_i$ and $r \in R$,
\begin{align*}
F_{\bigoplus M_i}(r(m_i)_{i \in I})
= (F_{M_i}(rm_i))_{i \in I}
= (d(r)m_i + rF_{M_i}(m_i))_{i \in I}
= d(r)(m_i)_{i \in I} + rF_{\bigoplus M_i}((m_i)_{i \in I}),
\end{align*}
and if $(m_i)_{i \in I}$ is homogeneous of degree $\tau$, each nonzero $m_i \in M_{i,\tau}$, so $F_{M_i}(m_i) \in M_{i,\tau}$ by hypothesis.

(2) Since $R$ is commutative, $M$ and $N$ are $(R,R)$-bimodules via $mr:=rm$, $nr:=rn$, so $M\otimes_R N$ is defined without ambiguity. The biadditive map $\Theta(m,n) := F_M(m)\otimes n + m\otimes F_N(n)$ satisfies, for $r\in R$,
\[
\Theta(mr,n) = F_M(mr)\otimes n + mr\otimes F_N(n)
= (d(r)m+rF_M(m))\otimes n + rm\otimes F_N(n)
= d(r)(m\otimes n) + r\,\Theta(m,n),
\]
and the same computation with $M$ and $N$ exchanged gives $\Theta(m,rn) = d(r)(m\otimes n) + r\,\Theta(m,n)$ as well. Since $mr\otimes n = m\otimes rn$ in $M\otimes_R N$, the two expressions agree, so $\Theta$ descends to a well-defined additive map $F_{M\otimes N}$ with $F_{M\otimes N}(m\otimes n) = \Theta(m,n)$. For $r \in R$, $m \in M_\tau$, $n \in N_\sigma$,
\[
F_{M \otimes N}(r(m \otimes n))
= (d(r)m + rF_M(m)) \otimes n + rm \otimes F_N(n)
= d(r)(m \otimes n) + rF_{M \otimes N}(m \otimes n),
\]
and degree preservation follows from $F_M(m) \in M_\tau$, $F_N(n) \in N_\sigma$, $m \otimes n \in (M \otimes_R N)_{\tau\sigma}$.
\end{proof}

\begin{definition}\label{HJK1EYY}
A graded $R$-module homomorphism $\phi: M \too N$ is a \defn{gr-generalized homomorphism} if $\phi \circ F_M = F_N \circ \phi$ for $(F_M,d)_{h,M} \in p\mathfrak{Der}^{gh}_\Gamma(R,M)$ and $(F_N,d)_{h,N} \in p\mathfrak{Der}^{gh}_\Gamma(R,N)$.
\end{definition}

\begin{example}\label{8II49KQ}
Let $\{M_i\}_{i \in I}$ be a finite family of graded $R$-modules with direct sum $M = \bigoplus_{i \in I} M_i$. If each $M_i$ admits a gr-generalized derivation $(F_{M_i},d)_{h,M_i}$ with the same associated derivation $d$, then the canonical projections $\pi_j: M \too M_j$ are gr-generalized homomorphisms.
\end{example}

\begin{proposition}\label{QHG9JRS}
Let $\phi: M \too N$ be a surjective graded $R$-module homomorphism between $\Gamma$-graded modules such that $\ker(\phi)$ is a gr-differential submodule of $M$. Then there is a well-defined $Z(R) \cap R_e$-linear map
\[
\phi_*: p\mathfrak{Der}^{gh}_\Gamma(R,M) \too p\mathfrak{Der}^{gh}_\Gamma(R,N)
\]
such that, for any $(F_M,d)_{h,M}$ with $F_M(\ker\phi) \subseteq \ker\phi$, the induced map is $(F_N,d)_{h,N} = \phi_*((F_M,d)_{h,M})$.
\end{proposition}

\begin{proof}
For $n \in N$, choose $m \in M$ with $\phi(m) = n$ and set $F_N(n) := \phi(F_M(m))$. If $\phi(m_1) = \phi(m_2)$, then $m_1-m_2\in\ker\phi$, so $F_M(m_1-m_2)\in\ker\phi$ by hypothesis, and $\phi(F_M(m_1)) = \phi(F_M(m_2))$: $F_N$ is well defined. For $r\in R$ and $m$ with $\phi(m)=n$,
\[
F_N(rn) = F_N(\phi(rm)) = \phi(F_M(rm)) = \phi(d(r)m + rF_M(m)) = d(r)n + rF_N(n).
\]
If $n \in N_\tau$, choose $m \in M_\tau$ with $\phi(m) = n$, using that $\phi$ is graded. Then $F_M(m) \in M_\tau$, so $F_N(n) = \phi(F_M(m)) \in N_\tau$. Linearity of $\phi_*$ over $Z(R) \cap R_e$ follows from the linearity of $\phi$.
\end{proof}

\begin{corollary}\label{NREK8PK}
For a graded isomorphism $\phi: M \too N$, the induced map
\[
\phi_*: p\mathfrak{Der}^{gh}_\Gamma(R,M) \too p\mathfrak{Der}^{gh}_\Gamma(R,N)
\]
is a $Z(R) \cap R_e$-module isomorphism, with inverse $\psi_*((F_N,d)_{h,N}) = (\phi^{-1} \circ F_N \circ \phi, d)_{h,M}$.
\end{corollary}

\begin{proof}
Since $\phi$ is a graded isomorphism, $\ker\phi=0$ is trivially a gr-differential submodule, so Proposition \ref{QHG9JRS} applies to $\phi$ and, symmetrically, to $\phi^{-1}$. The formula $F_N(\phi(m))=\phi(F_M(m))$ shows that $\psi_*(\phi_*((F_M,d)_{h,M}))$ agrees with $(F_M,d)_{h,M}$ on every $m\in M$, and symmetrically for $\phi_*\circ\psi_*$; hence $\psi_*$ is a two-sided inverse of $\phi_*$. Linearity of $\phi_*$ over $Z(R)\cap R_e$ is Proposition \ref{QHG9JRS}.
\end{proof}

\begin{definition}\label{3YR9OUD}
The category $\mathscr{M}^{gh}_\Gamma$ has:
\begin{enumerate}
\item \textit{Objects}: triples $(R, M, (F_M,d)_{h,M})$ where $R$ is $\Gamma$-graded, $M$ is a graded $R$-module, and $(F_M,d)_{h,M} \in p\mathfrak{Der}^{gh}_\Gamma(R,M)$;
\item \textit{Morphisms}: pairs $(\phi,\psi): (R, M, (F_M,d)_{h,M}) \too (S, N, (F_N,e)_{h,N})$, where $\phi: R \too S$ is a graded ring homomorphism, $\psi: M \too N$ is $\phi$-semilinear, and the diagrams
\[
\begin{tikzcd}
M \arrow[r, "\psi"] \arrow[d, "F_M"'] & N \arrow[d, "F_N"] \\
M \arrow[r, "\psi"'] & N
\end{tikzcd}
\quad\text{and}\quad
\begin{tikzcd}
R \arrow[r, "\phi"] \arrow[d, "d"'] & S \arrow[d, "e"] \\
R \arrow[r, "\phi"'] & S
\end{tikzcd}
\]
commute.
\end{enumerate}
\end{definition}

\begin{theorem}\label{P4E30MY}
$\mathscr{M}^{gh}_\Gamma$ is a well-defined category.
\end{theorem}

\begin{proof}
Let $(\phi,\psi): (R, M, (F_M,d)_{h,M}) \too (S, N, (F_N,e)_{h,N})$ and $(\phi',\psi'): (S, N, (F_N,e)_{h,N}) \too (T, P, (F_P,f)_{h,P})$ be morphisms in $\mathscr{M}^{gh}_\Gamma$. Then $\psi'\circ\psi$ is $(\phi'\circ\phi)$-semilinear:
\[
(\psi' \circ \psi)(rm) = \psi'(\phi(r)\psi(m)) = \phi'(\phi(r))\psi'(\psi(m)) = (\phi' \circ \phi)(r)(\psi' \circ \psi)(m),
\]
and compatibility with the derivations holds since
\[
(\psi' \circ \psi)\circ F_M = \psi'\circ(F_N\circ\psi) = F_P\circ(\psi'\circ\psi), \quad (\phi' \circ \phi)\circ d = f\circ(\phi'\circ\phi).
\]
So $(\phi'\circ\phi,\psi'\circ\psi)$ is again a morphism. For any object $(R, M, (F_M,d)_{h,M})$, the pair $(\text{id}_R, \text{id}_M)$ satisfies
\[
\text{id}_M(rm) = rm = \text{id}_R(r)\text{id}_M(m),\quad
\text{id}_M \circ F_M = F_M \circ \text{id}_M,\quad
\text{id}_R \circ d = d \circ \text{id}_R,
\]
and is therefore the identity morphism. Associativity and the identity laws follow directly from the corresponding properties of function composition.
\end{proof}

\printbibliography

\bigskip
\noindent
{
Sidi Mohamed Ben Abdellah University, Fez, Morocco\\
{\itshape e-mail:} \href{mailto:y.aitmohamed@yahoo.com}{y.aitmohamed@yahoo.com}

\end{document}